\documentclass[11pt]{amsart}

\usepackage{graphicx}

\newtheorem{theorem}{Theorem}[section]
\newtheorem{proposition}[theorem]{Proposition}
\newtheorem{lemma}[theorem]{Lemma}

\theoremstyle{definition}
\newtheorem{definition}[theorem]{Definition}

\numberwithin{equation}{section}

\newcommand{\depth}{\operatorname{depth}}

\renewcommand{\P}{\operatorname{{\mathcal P}}}

\newcommand{\D}{\operatorname{{\mathcal D}}}
\newcommand{\G}{\operatorname{{\mathcal G}}}

\newcommand{\T}{\operatorname{{\mathcal T}}}

\begin{document}

\title[Constructing knot tunnels using giant steps]{Constructing knot
tunnels using giant steps}

\author{Sangbum Cho}
\address{Department of Mathematics\\
University of Oklahoma\\
Norman, Oklahoma 73019\\
USA}
\email{scho@ou.edu}

\author{Darryl McCullough}
\address{Department of Mathematics\\
University of Oklahoma\\
Norman, Oklahoma 73019\\
USA}
\email{dmccullough@math.ou.edu}
\urladdr{www.math.ou.edu/$_{\widetilde{\phantom{n}}}$dmccullough/}
\thanks{The second author was supported in part by NSF grant
DMS-0802424}


\date{\today}


\begin{abstract}
In \cite{GST}, Goda, Scharlemann, and Thompson described a general
construction of all tunnels of tunnel number $1$ knots using ``tunnel
moves''. The theory of tunnels introduced in \cite{CM} provides a
combinatorial approach to understanding tunnel moves. We use it to
calculate the number of distinct minimal sequences of such moves that can
produce a given tunnel. As a consequence, we see that for a sparse infinite
set of tunnels, the minimal sequence is unique, but generically a tunnel
will have many such constructions.\par
\end{abstract}

\maketitle

\section*{Introduction}
\label{sec:intro}

In previous work \cite{CM}, we introduced a theory of tunnels of tunnel
number $1$ knots based on the disk complex of the genus-$2$ handlebody. It
provides a simplicial complex $\D(H)/\mathcal G$ whose vertices correspond
to the (equivalence classes of) tunnels of all tunnel number~$1$ knots.  As
we will explain below, two tunnels span a $1$-simplex of $\D(H)/\mathcal G$
exactly when each is obtained from the other by a construction given by
H. Goda, M. Scharlemann, and A. Thompson in~\cite{GST}. For reasons that
will become apparent, we call these constructions ``giant steps.'' The
connectivity of $\D(H)/\mathcal G$ shows that every tunnel can be obtained
from the unique tunnel $\pi_0$ of the trivial knot by some sequence of
giant steps, a fact already proved in~\cite{GST}.

In this note, we will use the combinatorial structure of $\D(H)/\mathcal G$
to examine minimal length sequences of giant steps that start from $\pi_0$
and produce a given tunnel $\tau$. Our main result is an algorithm to
calculate the number of distinct such sequences. In fact this number is
just the number of shortest paths between two vertices of the Farey
graph. The algorithm is quite elementary, but we have been unable to find
it in the literature. We will use it to see that for a sparse infinite set
of tunnels, the minimal construction sequence is unique, but generically a
tunnel will have many such sequences. The algorithm is effective and we
have implemented it computationally~\cite{slopes}.

The length of a minimal sequence of giant steps producing a given tunnel is
equal to the invariant called the \textit{depth} of the tunnel, defined
below. This invariant is used extensively in our work on bridge numbers of
tunnel number $1$ knots in~\cite{CMdepth}.

Only a minimal amount of the theory from \cite{CM} is needed for the
present application. We review it briefly in
Sections~\ref{sec:disk_complex}
and~\ref{sec:binary}. Section~\ref{sec:GST_moves} defines giant steps
precisely, and Section~\ref{sec:GST_move_seqs} presents the algorithm and
some of its uses.

\section{The tree of knot tunnels}
\label{sec:disk_complex}

\begin{figure}
\begin{center}
\includegraphics[width=4cm]{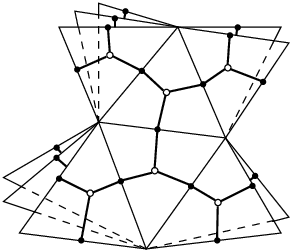}
\caption{A portion of the nonseparating disk complex $\D(H)$ and the tree
$\widetilde{\T}$. Countably many $2$-simplices meet along each edge.}
\label{fig:subdivision}
\end{center}
\end{figure}

Let $H$ be a genus~$2$ orientable handlebody, regarded as the standard
unknotted handlebody in $S^3$. For us, a \textit{disk in H} means a
properly imbedded disk in $H$, \textit{which is assumed to be nonseparating
unless otherwise stated.} The \textit{disk complex} $\D(H)$ is a
$2$-dimensional, contractible simplicial complex, whose vertices are the
isotopy classes of disks in $H$, such that a collection of $k+1$ vertices
spans a $k$-simplex if and only if they admit a set of pairwise-disjoint
representatives. Each $1$-simplex of $\D(H)$ is a face of countably many
$2$-simplices. As suggested by Figure~\ref{fig:subdivision}, $\D(H)$ grows
outward from any of its $2$-simplices in a treelike way. In fact, it
deformation retracts to the tree $\widetilde{\T}$ seen in
Figure~\ref{fig:subdivision}.

A tunnel of a tunnel number 1 knot produces a disk in $H$ as follows. The
tunnel is a $1$-handle attached to a regular neighborhood of the knot to
form an unknotted genus-$2$ handlebody. An isotopy moving this handlebody
to $H$ carries a cocore $2$-disk of that $1$-handle to a nonseparating disk
in $H$, and carries the tunnel number~$1$ knot to a core circle of the
solid torus obtained by cutting $H$ along that disk.

The indeterminacy in the choice of the isotopy is the group of isotopy
classes of orientation-preserving homeomorphisms of $S^3$ that preserve
$H$. This group is called the \textit{Goeritz group $\G$.} Work of
M. Scharlemann \cite{ScharlemannTree} and E. Akbas \cite{Akbas} proves that
$\G$ is finitely presented, and even provides a simple presentation of it.

\begin{figure}\begin{center}
\includegraphics[width=4cm]{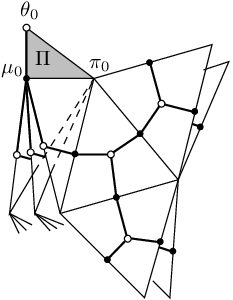}
\caption{A portion of $\D(H)/\G$ and $\T$ near the primitive orbits.}
\label{fig:Delta}
\end{center}
\end{figure}

Since two disks in $H$ determine equivalent tunnels exactly when they
differ by an isotopy moving $H$ through $S^3$, \textit{the collection of
all (equivalence classes of) tunnels of all tunnel number~$1$ knots
corresponds to the set of orbits of vertices of $\D(H)$ under $\G$.} So it
is natural to examine the quotient complex $\D(H)/\G$, which is illustrated
in Figure~\ref{fig:Delta}.

Through work of the first author~\cite{Cho}, the action of $\G$ on $\D(H)$
is well-understood. A \textit{primitive} disk in $H$ is a disk $D$ such
that there is a disk $E$ in $\overline{S^3-H}$ for which $\partial D$ and
$\partial E$ intersect transversely in one point in $\partial H$.  The
primitive disks (regarded as vertices) span a contractible subcomplex
$\P(H)$ of $\D(H)$, called the \textit{primitive subcomplex}. The action of
$\G$ on $\P(H)$ is as transitive as possible, indeed the quotient
$\P(H)/\G$ is a single $2$-simplex $\Pi$ which is the image of any
$2$-simplex of the first barycentric subdivision of $\P(H)$. Its vertices
are $\pi_0$, the orbit of all primitive disks, $\mu_0$, the orbit of all
pairs of disjoint primitive disks, and $\theta_0$, the orbit of all triples
of disjoint primitive disks. Regarded as a tunnel, $\pi_0$ is the unique
tunnel of the trivial knot.

On the remainder of $\D(H)$, the stabilizers of the action are as small as
possible. A $2$-simplex which has two primitive vertices and one
nonprimitive is identified with some other such simplices, then folded in
half and attached to $\Pi$ along the edge $\langle \mu_0,\pi_0\rangle$.
The nonprimitive vertices of such $2$-simplices are exactly the disks in
$\D(H)$ that are disjoint from some primitive pair, and these are called
$\textit{simple}$ disks. As tunnels, they are the upper and lower tunnels
of $2$-bridge knots, and we call them the simple tunnels. The remaining
$2$-simplices of $\D(H)$ receive no self-identifications, and descend to
portions of $\D(H)/\G$ that are treelike and are attached to one of the
edges $\langle \pi_0,\tau_0\rangle$ where $\tau_0$ is simple.

The tree $\widetilde{\T}$ shown in Figure~\ref{fig:subdivision} is
constructed as follows.  Let $\D'(H)$ be the first barycentric subdivision
of $\D(H)$. Denote by $\widetilde{\T}$ the subcomplex of $\D'(H)$ obtained
by removing the open stars of the vertices of $\D(H)$. It is a bipartite
graph, with ``white'' vertices of valence $3$ represented by triples and
``black'' vertices of (countably) infinite valence represented by
pairs. The valences reflect the fact that moving along an edge from a
triple to a pair corresponds to removing one of its three disks, while
moving from a pair to a triple corresponds to adding one of infinitely many
possible third disks to a pair.

\begin{figure}
\begin{center}
\includegraphics[width=5cm]{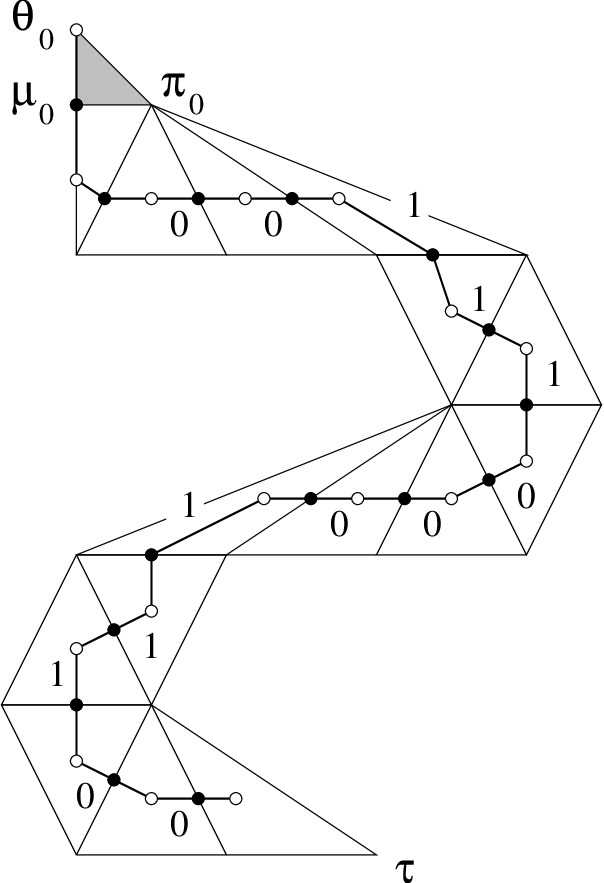}
\caption{The principal path of a tunnel $\tau$ having binary invariants
$0011100011100$, or equivalently with step sequence ``DRRRDRDLLLDLDRR''.}
\label{fig:ppath}
\end{center}
\end{figure}
The image $\widetilde{\T}/\G$ of $\widetilde{\T}$ in $\D'(H)/\G$ is a tree
$\T$. The vertices of $\D'(H)/\G$ that are images of vertices of $\D(H)$
are not in $\T$, but their links in $\D'(H)/\G$ are subcomplexes of
$\T$. These links are infinite trees. For each such vertex $\tau$ of
$\D'(H)/\G$, i.~e.~each tunnel, there is a unique shortest path in $\T$ from
$\theta_0$ to \textit{the vertex in the link of $\tau$ that is closest to
$\theta_0$.} This path is called the \textit{principal path} of $\tau$,
and this closest vertex is a triple, called the \textit{principal vertex}
of $\tau$. The two disks in the principal vertex, other than $\tau$, are
called the \textit{principal pair} of $\tau$. They are exactly the disks
called $\mu^+$ and $\mu^-$ that play a key role
in~\cite{Scharlemann-Thompson}. Figure~\ref{fig:ppath} shows the
principal path of a certain tunnel.

\section{The binary invariants}
\label{sec:binary}

Moving through the tree $\T/\G$ corresponds to a useful construction of
tunnels, called the cabling construction, but we will not need it here.
What is important for us is the combinatorial structure near the principal
path of $\tau$. This structure is determined by a sequence of ``binary''
invariants $s_2, s_3,\ldots\,$,~$s_n$, defined in \cite{CM}. We do not need
their formal definition, which involves the cabling construction, for we
can think of them in a very simple way, from the viewpoint of a traveler
along the path. A \textit{step} of the principal path is a portion between
successive white vertices. At Step~$0$, a traveler goes from $\theta_0$
to the principal vertex of some simple tunnel. At Step~$1$, the traveler
(whom we are viewing from above) must make a left turn. Starting with
Step~$2$, the traveler must make a choice of turning left or turning right
out of the white vertex.  The invariant $s_2$ is $0$ if this is a left turn
and $1$ if it is a right turn. In general, $s_i$ is $0$ if the direction of
the turn at Step~$i$ is the same as the direction of the turn at
Step~$i-1$, and $1$ if it is different.  Figure~\ref{fig:ppath} shows the
principal path of a tunnel with binary invariants $0011100011100$.\par

It is sometimes useful to describe the principal path from the reader's
viewpoint. The initial step is always down (``D'') and the second step, due
to the standard way that we draw the picture, is to the reader's right
(``R''). Each subsequent step is either left (``L''), down, or right.  An
``L'' can only be followed by another ``L'' or a ``D'', according as the
corresponding binary invariant $s$ is $0$ or $1$, and similarly an ``R'' is
followed by another ``R'' or a ``D'', according as $s$ is $0$ or $1$. When
the previous step is ``D'', then the effect of $s$ depends on the step
before the ``D''. If the two previous steps were ``LD'', then the next step
is ``R'' or ``L'' according as $s$ is $0$ or $1$, while if they were
``RD'', then the next step is ``L'' or ``R'' according as $s$ is $0$ or
$1$. For the example of Figure~\ref{fig:ppath}, the step sequence is
``DRRRDRDLLLDLDRR''.

There are simple algorithms for translating between these two descriptions,
and functions that do this are included in the software at~\cite{slopes}.

\section{Giant steps}
\label{sec:GST_moves}

\begin{definition} Let $\tau$ and $\tau'$ be tunnels.
We say that $\tau'$ is obtained from $\tau$ by a \textit{giant step} if
$\tau$ and $\tau'$ are the endpoints of a $1$-simplex of
$\D(H)/\G$. Equivalently, $\tau$ and $\tau'$ can be represented by disjoint
disks in~$H$.
\end{definition}

In \cite{GST}, Goda, Scharlemann, and Thompson gave a geometric definition
of giant steps (this is one reason for our selection of the name Giant
STeps), as follows. Let $\tau$ be a nonseparating disk in $H$, and let $K$
be a simple closed curve in $\partial H$ that intersects $\tau$
transversely in one point. Let $N$ be a regular neighborhood in $H$ of
$K\cup \tau$. Then the frontier of $N$ separates $H$ into two solid tori,
one a regular neighborhood of $K$, so $K$ is a tunnel number~$1$ knot.

In the previous construction, the meridian disk $\tau'$ of the solid torus
that does not contain $K\cup\tau$ is the \textit{unique} nonseparating disk
$\tau'$ in $H$ that is disjoint from $K\cup \tau$, and $\tau'$ is a tunnel
of $K$. That is, the construction produces a specific tunnel of the
resulting knot $K$. A giant step as we have defined it simply amounts to
choosing the $\tau'$ first; $K$ is then determined up to isotopy in $H$ and
in $S^3$, although not up to isotopy in~$\partial H$.

Since the complex $\D(H)/\G$ is connected, we have the following, which is
part of Proposition~1.11 of~\cite{GST}.
\begin{proposition}
Let $\tau$ be a tunnel of a tunnel number $1$ knot. Then there is a
sequence of giant steps that starts with the tunnel of the trivial knot and
ends with~$\tau$.\par
\label{prop:GST_sequence_exists}
\end{proposition}

The \textit{depth} of a tunnel $\tau$ is defined to be the distance in the
$1$-skeleton of $\D(H)/\G$ from $\pi_0$ to $\tau$. That is, the depth is
exactly the length of a minimal sequence of giant steps from $\pi_0$
to~$\tau$.

\section{Minimal sequences of giant steps}
\label{sec:GST_move_seqs}

In this section we give the algorithm to calculate the number of minimal
length sequences of giant steps that start from $\pi_0$, the tunnel of the
trivial knot, and end with a given tunnel~$\tau$. This is an elementary
combinatorial problem, and the reader will note that it is essentially the
problem of computing the number of distinct geodesics bewtween two points
in the Farey graph. We will use the algorithm to see that for a sparse
infinite set of tunnels, the minimal giant step sequence construction is
unique, but generically a tunnel will have many such constructions.

By a \textit{path} (between two vertices) in $\D(H)/\G$, we mean a
simplicial path in the $1$-skeleton of $\D(H)/\G$, passing through a
sequence of vertices that are images of vertices of $\D(H)$ (i.~e.~vertices
that represent tunnels).  We describe such a path simply by listing the
vertices through which it passes. From Section~\ref{sec:GST_moves}, we know
that the minimal sequences of giant steps from the $\pi_0$ to a given
tunnel $\tau$ correspond exactly to the minimal-length paths in $\D(H)/\G$
from $\pi_0$ to $\tau$. We will only be interested in minimal-length paths.
\begin{definition}
Let $\tau$ be a nontrivial tunnel. Define the \textit{corridor} of $\tau$,
$C(\tau)$, as follows. Write the vertices of the principal path of $\tau$
as $\theta_0$, $\mu_0$, $\mu_0\cup \tau_0$, $\mu_1\,$, $\mu_1\cup
\tau_1,\ldots\,$, $\mu_n\cup \tau_n$, where $\tau=\tau_n$. Then $C(\tau)$
is the union of the $2$-simplices whose barycenters are the $\mu_i\cup
\tau_i$ for $0\leq i\leq n$ (where $\mu_0\cup \tau_0$ is regarded as the
barycenter of the $2$-simplex spanned by $\pi_0$, $\mu_0$, and $\tau_0$).
\label{def:corridor}
\end{definition}
\begin{figure}
\begin{center}
\includegraphics[width=11cm]{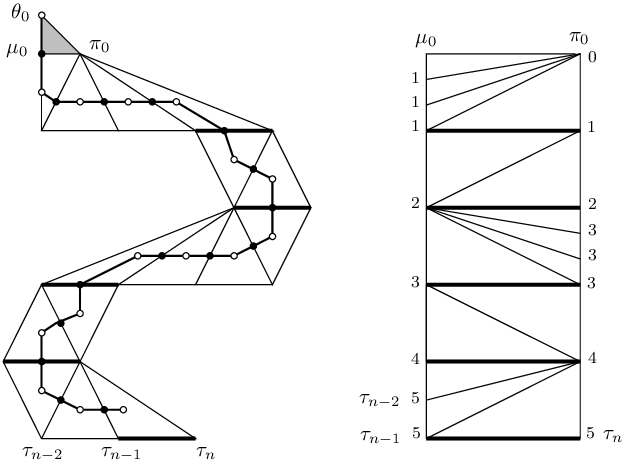}
\caption{The principal path and corridor $C(\tau_n)$ of a tunnel
$\tau_n$. The $\nabla$-edges are emphasized. In the picture of
$C(\tau_n)$ on the right, the depths of the tunnels are labeled.}
\label{fig:corridor}
\end{center}
\end{figure}

When $\tau$ is a simple tunnel, $C(\tau)$ is the triangle
$\langle\pi_0,\mu_0,\tau_0\rangle$. Otherwise, it can be viewed as a
rectangular or trapezoidal strip with top and bottom edges $\langle
\mu_0,\pi_0\rangle$ and $\langle \tau_{n-1},\tau_n\rangle$, as in
the drawing on the right in Figure~\ref{fig:corridor}.
\begin{lemma} Let $\tau$ be a tunnel, and let
$\sigma_0$, $\sigma_1,\ldots\,$, $\sigma_n$ be a path in $\D(H)/\G$ of
minimal length among the paths connecting the vertices $\sigma_0$ to
$\sigma_n$. If $\sigma_0$ and $\sigma_n$ lie in $C(\tau)$, then each
$\sigma_i$ lies in~$C(\tau)$.\par
\label{lem:stay_in_corridor}
\end{lemma}

\begin{proof}
If the lemma is false, then there exist $i$ and $j$ with $0\leq i<i+1<j\leq
n$ for which $\sigma_i$ and $\sigma_j$ lie in $C(\tau)$, but $\sigma_k$
does not lie in $C(\tau)$ for any $k$ with $i<k<j$.

The vertex $\{ \sigma_i,\sigma_{i+1}\}$ lies in the link $L$ of
$\sigma_i$ in $\D'(H)/\G$. This link is a tree, so there exists a vertex
$\sigma_i'$ in $C(\tau)$ such that $\{ \sigma_i,\sigma_i'\}$ is
connected to $\{ \sigma_i,\sigma_{i+1}\}$ by a path in $L$
meeting $C(\tau)$ only in $\{\sigma_i,\sigma_i'\}$. The
$1$-simplex in $\D(H)/\G$ spanned by $\sigma_i$ and $\sigma_i'$ separates
$\D(H)/\G$, with $C(\tau)$ and $\sigma_{i+1}$ lying in different
components. Therefore $\sigma_j$ must equal either $\sigma_i$ or
$\sigma_i'$. In either case we obtain a shorter path from $\sigma_0$ to
$\sigma_n$.
\end{proof}

In the special case that $\tau$ is of depth $1$, $\tau$ lies in the link in
$\D(H)/\G$ of $\pi_0$, and there is a unique path of length~$1$ from
$\pi_0$ to $\tau$. From now on, we assume that $\tau$ has depth at
least~$2$.

Now, regard $C(\tau)$ as in the diagram on the right in
Figure~\ref{fig:corridor}, with the edge $\langle \mu_0, \pi_0\rangle$ on
top, and with $\tau$ as one of the endpoints of the bottom edge.  In the
triangulation of $C(\tau)$, a \textit{$\nabla$-edge of depth $i$} is an
edge whose endpoints have depth $i$ and lie on different sides of
$C(\tau)$, and for which all vertices lying below its endpoints on either
side have depth greater than $i$. In Figure~\ref{fig:corridor}, the
$\nabla$-edges are highlighted.

Since the endpoints of any edge of $C(\tau)$ can have depths that differ by
at most $1$, there exists a unique $\nabla$-edge $\nabla(i)$ in $C(\tau)$
of depth $i$ for each $i$ with $1\leq i<\depth(\tau)$. There is also a
$\nabla$-edge for $i=\depth(\tau)$, unless $\tau$ and the endpoints of
$\nabla(\depth(\tau)-1)$ span a $2$-simplex.

The name $\nabla$-edge arises from the fact that (except for
$\nabla(\depth(\tau))$) the $\nabla$-edges are the tops of $2$-simplices of
the corridor that appear as $\nabla$'s when the corridor is drawn with
depth corresponding to the vertical coordinate, as in the diagram on the
left in Figure~\ref{fig:corridor}. Every nonprimitive $2$-simplex of
$\D(H)/\G$ has two vertices of the same depth and a third of depth either
larger by $1$ or smaller by $1$ than that common depth; for a ``$\nabla$''
$2$-simplex that depth is larger by $1$, while it is smaller by $1$ for a
``$\Delta$'' $2$-simplex.

Denote the left and right endpoints of $\nabla(i)$ by
$\partial_L(\nabla(i))$ and $\partial_R(\nabla(i))$ respectively.
\begin{lemma} Let $\nabla(i-1)$ and $\nabla(i)$ be successive
$\nabla$-edges. Then at least one of the pairs
$\{\partial_L(\nabla(i-1)),\partial_L(\nabla(i))\}$ and
$\{\partial_R(\nabla(i-1)),\partial_R(\nabla(i))\}$ are
the endpoints of an edge that lies in a side of~$C(\tau)$.
\label{lem:side_edge}
\end{lemma}
\begin{proof}
For each endpoint of $\nabla(i)$, select a path of length $i$ from the
endpoint to $\pi_0$.  By Lemma~\ref{lem:stay_in_corridor}, these paths lie
in $C(\tau)$. In particular, each of their first edges connects an endpoint of
$\nabla(i)$ to an endpoint of $\nabla(i-1)$. At most one of these first
edges can be diagonal, so at least one lies in a side.
\end{proof}

\begin{figure}
\begin{center}
\includegraphics{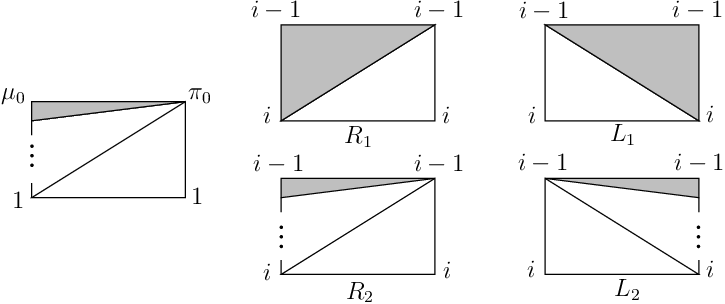}
\caption{The configuration above $\nabla(1)$, and the four possible
configurations between two $\nabla$-edges. In $R_1$ and $L_1$ there is
only one $2$-simplex above the diagonal edge, while in $R_2$ and $L_2$
there are two or more. In the configuration above $\nabla(1)$, there may be
only one $2$-simplex above the diagonal. The shaded $2$-simplices are
$\nabla$ $2$-simplices. The letter $L$ (respectively, $R$) signifies that
the portion below the $\nabla$ simplex contains ``L'' steps
(respectively, ``R'' steps) of the principal path.}
\label{fig:configurations}
\end{center}
\end{figure}

Lemma~\ref{lem:side_edge} shows that the triangulation of the portion of
$C$ between $\nabla(i-1)$ and $\nabla(i)$ must have one of the four
configurations $L_1$, $R_1$, $L_2$, or $R_2$ shown in
Figure~\ref{fig:configurations}. The portion of $C(\tau)$ above $\nabla(1)$
must be as in the leftmost diagram in Figure~\ref{fig:configurations},
where there may be only one $2$-simplex above the diagonal.

Now, we show how to calculate the number of minimal paths from $\pi_0$ to
$\tau$. Denote by $\lambda_i$ the number of paths in $C(\tau)$ of length
$i$ from $\pi_0$ to the left endpoint of $\nabla(i)$, and by $\rho_i$ the
number to its right endpoint. Clearly $\begin{pmatrix} \lambda_1 \\ \rho_1
\end{pmatrix}=\begin{pmatrix} 1 \\ 1 \end{pmatrix}$.

Let $k$ be the maximum $i$ for which $\nabla(i)$ is defined.  For each
$2\leq i\leq k$, let $C_i$ be $L_1$, $R_1$, $L_2$, or $R_2$ according to
which of the four configurations in Figure~\ref{fig:configurations}
describes the triangulation of $C(\tau)$ between $\nabla(i-1)$ and
$\nabla(i)$. For $2\leq i\leq k$, put $M_i$ equal to the matrix given in
the following table, according to the value of $C_i$:
\medskip

\begin{small}
\setlength{\fboxsep}{0pt}
\setlength{\fboxrule}{0.5pt}
\begin{center}
\fbox{%
\begin{tabular}{c|c|c|c|c}
$C_i$&$L_1$&$R_1$&$L_2$&$R_2$\\
\hline
$M_i$&$\begin{pmatrix}1 & 0\\ 1 & 1 \end{pmatrix}$
&$\begin{pmatrix}1 & 1\\ 0 & 1 \end{pmatrix}$
&$\begin{pmatrix}1 & 0\\ 1 & 0 \end{pmatrix}$
&$\begin{pmatrix}0 & 1\\ 0 & 1 \end{pmatrix}$\\
\end{tabular}}
\end{center}
\end{small}
\medskip

Observe that
\[ M_i\begin{pmatrix} \lambda_{i-1} \\ \rho_{i-1}\end{pmatrix} =
\begin{pmatrix} \lambda_i \\ \rho_i\end{pmatrix}\ .\]
Therefore we have
\[ \begin{pmatrix} \lambda_k \\ \rho_k\end{pmatrix} =
M_kM_{k-1}\cdots M_2\begin{pmatrix} 1 \\ 1 \end{pmatrix}\ .\]

If $\tau$ and the endpoints of $\nabla(k)$ span a $2$-simplex, as in the
case of the tunnel $\tau_{n-2}$ in Figure~\ref{fig:corridor}, then the
number of minimal-length paths from $\pi_0$ to $\tau$ is
$\lambda_k+\rho_k$.  Otherwise, $\tau$ is the left or right endpoint of
$\nabla(k)$, and the number of paths equals $\lambda_k$ or $\rho_k$
respectively.

The algorithm just described is not difficult to implement
computationally~\cite{slopes}. For the example in
Figure~\ref{fig:corridor}, the output of the program is:
\medskip

\begin{ttfamily}
Depth> gst( '0011100011100', verbose = True )

The block configurations are R1, L2, L1, R2.

The transformation matrices M\_2, ..., M\_k are:

\ \ \ [ [ 1, 1 ], [ 0, 1 ] ]

\ \ \ [ [ 1, 0 ], [ 1, 0 ] ]

\ \ \ [ [ 1, 0 ], [ 1, 1 ] ]

\ \ \ [ [ 0, 1 ], [ 0, 1 ] ]

and M\_k * ... * M\_2 is [ [ 2, 2 ], [ 2, 2 ] ].

This tunnel has 4 minimal giant step constructions.
\end{ttfamily}
\medskip

Some examples are the tunnels whose parameter sequences are the following:
\begin{enumerate}
\item $s_2s_3\cdots s_n=100100\cdots 100$, whose corresponding step
sequence is $DRDLLDRRDLL \cdots DRR$ (or $\cdots DLL$). The configuration
sequence alternates as $L_2$, $R_2$, $L_2$, $R_2\ldots\,$, and there is a
unique minimal giant step sequence.
\item $s_2s_3\cdots s_{2n+1}=1010\cdots 10$, or $DRDLDRDLD\cdots DR$ (or
$\cdots DL$). The configuration sequence alternates as $L_1$, $R_1$, $L_1$,
$R_1\ldots\,$, and the number of minimal giant step sequences is the term
$F_n$ of the Fibonacci sequence $(F_0,F_1,F_2,\ldots)=(1,1,2,3,5,\ldots)$.
\item $s_2s_3\cdots s_{2n+1}=111\cdots 1$, an even number of $1$'s.  The
step sequence is $DRDRDR\cdots DR$. The configuration sequence is $R_1$,
$R_1,\ldots\,$, $R_1$, $\tau$ is the right-hand endpoint of $\nabla(n)$,
and there is a unique minimal giant step sequence.
\item $s_2s_3\cdots s_{2n}=111\cdots 1$, an odd number of $1$'s. The step
sequence is $DRDRDR\cdots D$. The configuration sequence is again $R_1$,
$R_1,\ldots\,$, $R_1$, but $\tau$ lies in a $\nabla$ $2$-simplex below
$\nabla(n)$, and there are $n+1$ minimal giant step sequences.
\end{enumerate}
Examples of the last two types are obtained from each other by a single
additional cabling construction, even though the numbers of minimal giant
step constructions differ by arbitrarily large amounts.

The algorithm shows that tunnels with a unique minimal giant step sequence
are sparse. For instance, the product (in reverse order) of the matrices
determined by the configuration sequence $L_1R_1L_1R_1$ has all entries
greater than~$1$, and whenever this product appears as any block of four
terms in the product $M_k* \cdots * M_2$ that occurs in the algorithm,
there must be more than one minimal giant step sequence. Configuration
sequences containing $L_1R_1L_1R_1$ are generic in any reasonable sense.

\bibliographystyle{amsplain}

\end{document}